\documentstyle[12pt,amstex,amssymb]{article}
\newcommand\qed{{\unskip\nobreak\hfil\penalty50\hskip2em\vadjust{}
    \nobreak\hfil$\Box$\parfillskip=0pt\finalhyphendemerits=0\par}}

\newtheorem{lemma}{Lemma}

\newcommand{\x}{\times}
\newcommand{\<}{\langle}
\renewcommand{\>}{\rangle}

\renewcommand{\a}{\alpha}

\newcommand{\G}{\Gamma}
\renewcommand{\l}{\lambda}
\renewcommand{\L}{\Lambda}

\newcommand{\s}{\sigma}

\newcommand{\th}{\theta}

\renewcommand{\i}{\infty}
\newcommand{\p}{\partial}

\title{A Topological Method to Compute Spectral Flow}
\author{Dave Auckly\thanks{Partially supported by an NSF Postdoctoral
Fellowship while the author was visiting MSRI. } }

\begin{document}

\maketitle

\bigskip\medskip
\begin{quote}
{\sc Abstract} \newline
This paper describes a topological method to compute the spectral flow of
a family of twisted Dirac operators, it includes two detailed examples. Briefly, a formula of Atiyah, Patodi and Singer expresses the spectral flow in terms of Chern-Simons invariants and rho invariants. The first step is to construct a
flat cobordism to a new bigger 3-manifold.  The advantage of the new connection
is that it is in the path component of a reducible connection. The second step
is to calculate the effect of these operations on the invariants.  The final
step is an application of the G-signature theorem to compute the invariants.
\end{quote}



The spectral flow of a family of operators is a generalization of the signature
of a self-adjoint finite dimensional operator.  If $A_t$ is a family
of operators with a real, discrete, spectrum, then the spectral flow of $A_t$ is
the number of eigenvalues that move from negative to positive minus the number
which move from positive to negative.  Recalling that the number of positive
eigenvalues minus the number of negative eigenvalues, we see that
\[
SF(A_t) = \mbox{ Sign } A_1 - \mbox{ Sign } A_0
\]
for finite dimensional operators. 

In this paper, we extend to spectral flow the method that we used to compute
Chern--Simons invariants of flat $SU_2$ connections in a previous paper
[A].  After outlining the method, we work out two examples in detail.
Briefly, a formula of Atiyah, Patodi and Singer expresses the spectral flow in
terms of Chern--Simons invariants and rho invariants [APS].  To compute the rho
invariant of a flat $SU_2$ connection on one $3$-manifold, we construct a flat
cobordant connection on a new, bigger $3$-manifold.  It is straightforward to
compute the difference between the rho invariant of a connection and a flat
cobordant connection.  The advantage of the new connection is that it is in the
path component of a reducible connection in the space of flat connections.  It
is possible to compute the change in the rho invariant along a path of flat
connections.  Finally, the reducible connection may be assumed to have finite
holonomy so its rho invariant may be computed with the $G$-signature theorem.

We begin by describing a collection of self-adjoint elliptic operators over any
$3$-manifold and stating the Atiyah--Patodi--Singer formula for the spectral
flow.  Let $A$ be a special unitary connection on a vector bundle over a
$3$-manifold, $M$.

\medskip
\noindent
{\bf Definition}.  $B(A): \G(\L^0 \oplus \L^2M \otimes {\Bbb C}^n) \to
\G(\L^0 \oplus \L^2M \otimes {\Bbb C}^n)$;
\[
\a \mapsto i^{\mbox{deg } \a}(*d_A\a - d_A * \a)
\]
is {\em the Chern--Simons Hessian}.

\medskip
The unitary connections that we consider will be constructed in the same way.
Start with a homomorphism from $\pi_1(M)$ to $SU_2$.  This induces a
representation into $U_3$ by composition with the map:
\[
\begin{array}{l}
\mbox{Ad}: SU_2 \to U_3; \\
(\mbox{Ad } A)(X \otimes z) = A X A^{-1} \otimes z.
\end{array}
\]
Here we are identifying ${\Bbb C}^3$ with $su_2 \otimes {\Bbb C}$.  Finally,
the canonical flat connection on the universal covering of $M$ induces a
natural connection on the twisted product, ${\tilde M} \x_{\pi_1(M)} U_3$.
Given two $SU_2$ representations, $\a_0$ and $\a_1$ we construct the
corresponding $U_3$ connections, $A_0$ and $A_1$, and the path of connections
\[
A_t = (1-t)A_0 + tA_1.
\]
By the spectral flow from $\a_0$ to $\a_1$, we will mean the spectral flow of
the family of operators $B(A_t)$.
When one of the operators $A_0$ or $A_1$ has
a non-trivial kernel, we need to pick some convention to count the zero modes.
We count the zero modes of $A_0$ as positive eigenvalues and the zero modes of
$A_1$ as negative eigenvalues.  This is the convention used by Fintushel and
Stern [FS1], [FS2].  After the first example, we will explain
why this is a natural convention.

As an aside, we will explain why the spectrum of $B(A_t)$ is discrete.  For any
constant, $\l$, $B(A_t) - \l$ is an elliptic operator.  It follows from
elliptic regularity and Rellich's Lemma that $B(A_t) - \l$ is Fredholm.  This
in turn implies that the spectrum of $B(A_t) - \l$ does not accumulate near
zero.

More than this is true.  Atiyah, Patodi and Singer show that the function
\[
\eta_{B(A)}(s) = \sum_{\l \in \mbox{ Spec}(B(A))} \l^{-s} \mbox{ Sign }
\l
\]
is an analytic function of $s$ for $\mbox{Re } s \gg 1$.  By analytic
continuation, they show that $\eta_{B(A)}(0)$ may be defined as a finite
number [APS].  This eta invariant is an extension of the signature to a
large class of infinite rank operators.  Let $\th$ denote the connection
constructed from the trivial representation.

\medskip
\noindent
{\bf Definition}.  The {\em rho invariant} of a representation, $\a$, is
\[
\rho_M(\a) = \eta_{B(A)} - \eta_{B(\th)}
\]
where $A$ is the connection constructed from $\a$.  The {\em Chern--Simons}
invariant of a representation is
\[
CS_M(\a) = \int_W - \frac {1}{4\pi^2} \mbox{ Tr}(\L^2F_{\bar A})
\]
where $W$ is any $4$-manifold with boundary $M$, and ${\bar A}$ is any $SU_2$
connection which restricts to the $SU_2$ connection induced from the
representation on the boundary.

The main theorem of Atiyah, Patodi and Singer may be stated in two different
ways depending on the boundary conditions.  If $D: \G(E) \to \G(F)$ is an
elliptic first order operator that restricts to the form
\[
D = \frac {\p}{\p t} + B
\]
on a metric product collar of the boundary of the base manifold, then one way
to put boundary conditions on the operator is to restrict the domain to
sections of $E$ that have $Ps|_{\p X} = 0$ where $P$ is the spectral projection
of $B$ corresponding to non-negative eigenvalues.  An alternate way to put
boundary conditions on $D$ is to introduce a non-compact base manifold, ${\hat
X} = X \cup_{\p X} ([0,\i) \x \p X)$ with a product metric on the end, and
first extend $D$ to an operator over ${\hat X}$, say ${\hat D}$, then restrict
the domain and range to $L^2$-sections.  In either case Atiyah, Patodi and
Singer prove that the index of $D$ is an integral of curvature terms plus a
correction term depending on the dimension of the kernel of $B$ and the eta
invariant of $B$ [APS].  We will call the first index the APS index
and the second the $L^2$ index.  A path of connections on $M$ describes a
connection, ${\bar A}$, on $[0,t] \x M$.  The APS theorem gives us the
following formula for the APS index of the signature operator twisted by ${\bar
A}$:
\[
\mbox{Index}(d_{\bar A} + d_{\bar A}^*|_{\Lambda^+}) = 16(CS_M(A_0) -
CS_M(A_t)) - h(A_0) - h(A_t) + \rho_M(A_0) - \rho_M(A_t).
\]
Here $h(A) = \dim \ker B(A)$.  Assuming that there is a consistent way to split
the spectrum of $B(A_t)$ into a finite part and an infinite part.  The spectral
flow of $B(A_t)$ may be expressed in terms of the finite part of the eta
invariants and the $h(B(A_t))$.  The above index formula can then be arranged
to set an integer valued function of the spectral flow and the $L^2$ index and
$h(B(A_t))$ equal to a continuous function of $t$.  Evaluating the continuous
function at $t = 0$ gives zero.  Finally, solving for the spectral flow gives:

\begin{lemma}
\label{lem1}
\[
SF(\a_0,\a_1) = 8(CS_M(\a_1) - CS_M(\a_0)) + \frac {1}{2} (\rho_M(\a_1) -
\rho_M(\a_0) - h(\a_1) - h(\a_0)).
\]
\end{lemma}

The general case of this result follows by cutting the path of connections into
parts where the spectrum may be split into a finite part and an infinite part.
This argument may be found in a paper by Kirk, Klassen and Ruberman
[KKR].

The Chern--Simon terms in the above formula may be computed by the method in
[A].  The connection corresponding to a representation is flat, so the
sequence
\[
\G(\L^0M \otimes {\Bbb C}^3) \stackrel{d_A}{\rightarrow} \G(\L^1M \otimes {\Bbb
C}^3) \stackrel{d_A}{\rightarrow} \G(\L^2M \otimes {\Bbb C}^3)
\stackrel{d_A}{\rightarrow} \G(\L^3M \otimes {\Bbb C}^3)
\]
is a complex.  Hodge theory implies that $h(\a) = \dim H^0(M;\a) \oplus
H^2(M;\a)$.  By Poincar\'e duality we see that $h(\a) = \dim H^0 \oplus H^1$.
Since the $2$-skeleton of $M$ is the $2$-skeleton of a $K(\pi_1(M),1)$, we see
that computing $h(\a)$ is reduced to an algebraic computation in group
cohomology.  In fact, $H^0(\pi_1(M);SU_2) = SU_2^{\pi_1(M)}$, thus
\[
\dim H^0 = \left\{ \begin{array}{ll}
0 &\mbox{if $\a$ is irreducible} \\
1 &\mbox{if $\a$ is abelian but not central} \\
3 &\mbox{if $\a$ is central.}
\end{array} \right.
\]
and $\dim H^1 = \dim Z^1 - 3 + \dim H^0$ where
\[
Z^1 = \{f: \pi_1(M) \to su_2 \mid f(gh) = f(g) + \a(g)f(h)\a(g^{-1})\}
\]
is the space of cocycles.  This leaves the rho invariants which may be computed
with a refined version of the same three steps which are used to compute
Chern--Simons invariants.

The first step in computing the rho invariant is to construct a flat cobordism
to a new representation and compute the change in the rho invariant.  The
standard cobordism to use is $W^4 = [0,1] \x M \cup_{S^1 \x D^2} V$, where $V =
[0,1] \x S^1 \x D^2/\sim$ $(0,\l,z) \qquad \sim (1,{\bar \l},{\bar z})$ and $S^1 \x
D^2 \to S^1 \x [1/3,2/3]^2 \to V$; \  $\l,s,t \mapsto (s,\l,\mbox{{\em exp}}(2\pi it))$.  The
difference between a pair of $L^2$ indices gives the following formula.

\begin{lemma}
\label{lem2}
If $\a: \pi_1(W^4) \to SU_2$, then
\[
\rho_{\p W}(\a|_{\p W}) = 3 \mbox{ {\em Sign} } W - \mbox{ {\em Sign} } Q_{\a}
\]
where $\mbox{{\em Sign }} W$ is the signature of the intersection form on
$\ker(H^2(W) \to H^2(\p W))$ and $\mbox{{\em Sign} } Q_{\a}$ is similar but
with twisted coefficients.
\end{lemma}

The group $\ker(H^2(W;\a) \to H^2(\p W;\a))$ is called the $L^2$-cohomology of
$W$, $H_{L^2}^2(W)$.  The following lemma computes the $L^2$ cohomology of the
standard cobordism.

\begin{lemma}
\label{lem3}
$H_{L^2}^2(W;\a)$ injects in \mbox{coker } $(H^1(M;\a) \oplus H^1(V,\a) \to H^1(S^1 \x
D^2;\a))$.
\end{lemma}

\medskip
\noindent
{\bf Proof}.  In the diagram,
\[
\begin{array}{rcl}
H^1([0,1] \x M) \oplus H^1(V) \to H^1(S^1 \x D^2) \to H^2(W) &\to &H^2([0,1] \x
M) \oplus H^2(V) \\
\downarrow &''' &\downarrow \\
H^2(\p W) &\to &H^2(M),
\end{array}
\]
the top row is exact and the map on the right is an injection. \qed

\medskip
The next step is to construct a path of representations from the new
representation to a reducible representation and compute the change in the rho
invariant across this path.  There is a general method to compute the spectral
flow along a path of flat connections due to Kirk and Klassen [KK2].  We
do not, however, need this formula, because we may pick a path with special
properties.

\begin{lemma}
\label{lem4}
Let $\a_t$ be a family of representations so that,
\[
\dim H^0(M;\a_t) = \left\{ \begin{array}{rl}
0 &\mbox{if $t > 0$} \\
1 &\mbox{if $t = 0$}
\end{array} \right.
\]
and,
\[
\dim H^1(M;\a_t) = \left\{ \begin{array}{rl}
h &\mbox{if $t > 0$} \\
h+1 &\mbox{if $t = 0$,}
\end{array} \right.
\]
then $\rho_M(\a_1) = \rho_M(\a_0)$.
\end{lemma}

\medskip
\noindent
{\bf Proof}.  We will use the formula in Lemma~1.  First compute,
\[
\begin{array}{rll}
CS(\a_1) - CS(\a_0) &= &-\frac {1}{4\pi^2} \int_{[0,1] \x M}
\mbox{ Tr}(\Lambda^2F_{\bar A}) \\
&= &-\frac {1}{4\pi^2} \int_{[0,1] \x M} \mbox{ Tr}(\Lambda^2(F_{\a_t} + dt
\wedge \frac {\p \a_t}{\p t})) \\
&= &-\frac {1}{4\pi^2} \int_{[0,1] \x M} \mbox{ Tr}(\Lambda^2dt \wedge \frac {\p
\a_t}{\p t}) = 0.
\end{array}
\]
Now, write
\[
\G((\L^0 \oplus \L^2)M \otimes {\Bbb C}^3) \cong \G(\L^0M \otimes {\Bbb C}^3)
\oplus \ker(d_A*) \oplus \mbox{ Im}(d_A*).
\]
We can do this because, $d_A*$ is self-adjoint.  This decomposes $B(A)$ as:
\[
B(A) = \left[ \begin{array}{ccc}
0 & -*d_A & 0 \\
*d_A & 0 & 0 \\
0 & 0 & d_A*
\end{array} \right] .
\]
It follows that the two eigenvalues that approach zero as $t \to 0$ have
opposite signs.  Just check that when $\left[ \begin{array}{c}
u \\
v \\
0
\end{array} \right]$ is an eigenvector with eigenvalue $\l$, $\left[
\begin{array}{c}
u \\
-v \\
0
\end{array} \right]$ is an eigenvector with eigenvalue $-\l$.  From this it is
apparent that $SF(\a_0,\a_1) = -h - 1$.  Plugging into Lemma~1 will now give
the result. \qed

\medskip
It is sometimes useful to compute the change in rho invariants along a path of
reducible connections.  In the typical case, the dimension of $\ker B(A_t)$
will be a constant, the spectral flow will be determined by the zero modes, and
the rho invariant will be constant.

The final step in the computation of the rho invariant and thus the spectral
flow, is to use the $G$-signature theorem to compute the rho invariant of the
reducible representation.  If $G$ is a finite group which acts effectively on a
$4$-manifold $W$, and $g \in G$, then we can define a $g$-signature.

\medskip
\noindent
{\bf Definition}.  $\mbox{Sign}(g,W) = \mbox{ Tr}(g^*|_{H_+^2(W)}) -
\mbox{ Tr}(g^*|_{H_-^2(W)})$.

\medskip
A group element may have a surface of fixed points modeled on, $g \cdot (z,w) =
(z,e^{i\psi}w)$ for $(z,w) \in {\Bbb C}^2$ or have isolated fixed points
modeled on $g \cdot (z,w) = (e^{i\th_1}z,e^{i\th_2}w)$.  Define a local
$g$-signature by
\[
L(g,W) = \sum_{F,p} F \cdot F \csc^2(\psi(F)/2) - \cot(\frac {\th_1(p)}{2})
\cot(\frac {\th_2(p)}{2}).
\]
The $g$-signature theorem states that,
\[
\mbox{Sign}(g,W) = L(g,W).
\]
If $G$ acts freely on a $3$-manifold, $M$, then there is a $4$-manifold with
boundary several copies of $M$ which the action extends over.  We can see this
by using the transfer map to prove that ${\bar H}_*(BG)$ is torsion for finite
groups.  In fact, by a direct geometric construction, we will see that we may
take $\p W = M$, when $G$ is finite cyclic.  In this setting define a
signature defect by
\[
\s_M(g) = L(g,W) - \mbox{ Sign}(g,W).
\]
This is independent of $W$ by the $g$-signature theorem.  If a representation,
factors through a finite group,
\[
\begin{array}{rcl}
\a:\pi_1(M) &\rightarrow &SU_2 \\
\searrow & &\nearrow \\
&G
\end{array}
\]
Then there is a formula for the rho invariant in terms of the signature
defect.  By an argument from [A], we may assume that every reducible
representation has this form.

\begin{lemma}
\label{lem5}
\[
\rho_M(\a) = \frac {1}{|G|} \sum_{g \ne 1} \s_{\hat M}(g)(\mbox{Tr } \mbox{Ad }
\a(g) - 3).
\]
(${\hat M}$ is the cover of $M$ coming from $\pi_1(M) \to G$.)
\end{lemma}

The proof of this may be found in [APS].  There is also an expression
for the signature defects involving the rho invariants.  The rho invariant may
be considered to be a generalization of an invariant for finite covering spaces
to infinite covering spaces.  The main idea of this paper is to go backward and
reduce the computation of this invariant to the case of finite coverings.  This
method works for any representation that is flat cobordant to a representation
in the path component of a reducible representation.  It works the best when
the path may be chosen to have the property in Lemma~3.  Every representation
on every graph manifold has this property.  Also some representations on some
hyperbolic manifolds have this property.  It is possible that every
representation on every $3$-manifold has this property.

Now that we have described a method to compute the spectral flow, we will work
out two examples.  The first example is a representation on the Poincar\'e
homology sphere, $\Sigma(2,3,5)$.  The second example is a representation on a
hyperbolic manifold.

The spectral flow for $\Sigma(2,3,5)$ could be computed with Lemma~1 and
Lemma~4.  It could also be computed with the Fintushel--Stern method by
considering the mapping cylinder down to the base orbifold.  Even so, it is a
simple example that well illustrates all of the important aspects of our method.

Since the fundamental group of $\Sigma(2,3,5)$ is fairly simple,
\[
\pi_1(\Sigma(2,3,5)) = \< Q_1,Q_2,Q_3,H \mid
[H,Q_i],Q_1^2H,Q_2^3H,Q_3^5H,HQ_1Q_2Q_3\>,
\]
it is not hard to find all $SU_2$ representations.  There are three
representations; the trivial one, and two irreducible representations.  The
fundamental representation is given by
\[
H \mapsto -1, Q_1 \mapsto i, Q_2 \mapsto g \exp(i\pi/3) g^{-1}, Q_3 \mapsto
h\exp(i\pi/5) h^{-1},
\]
where $g$ and $h$ are elements chosen to solve the final relation.  Here, we
are identifying $SU_2$ with $Sp_1$ and writing elements as unit quaternions.

We will compute the spectral flow from the trivial representation to the
fundamental representation.  The main part of the computation is computing the
rho invariant. The first step in the computation of the rho invariant is to
compute the change across the cobordism.

Since $\Sigma(2,3,5)$ is a Seifert fiber space, we can see the cobordism at the
level of base orbifolds.

\begin{figure}[h]
\vspace{2in}
\caption{The cobordism.} \label{fig1}
\end{figure}

\noindent
The base orbifold of the cobordism is a $2$-sphere with three cone points
cross with an interval glued along a disk to a solid Klein bottle.  A
handlebody decomposition of the part of the cobordism over the solid Klein
bottle is drawn in Figure~2.

\clearpage

\begin{figure}[h]
\vspace{2in}
\caption{The standard part of the cobordism.} \label{fig2}
\end{figure}

A neighborhood of the dotted curve above is identified with the regular fiber
in $\Sigma(2,3,5) \x \{1\} \subseteq \Sigma(2,3,5) \x [0,1]$.  This regular
fiber is drawn as a dotted circle in Figure~3.

\begin{figure}[h]
\vspace{2in}
\caption{A regular fiber.} \label{fig3}
\end{figure}

By Lefschitz duality, we may compute twisted signatures in either
$L^2$-cohomology or $L^2$-homology.  By Lemma~\ref{lem3}, the elements of the
second
$L^2$-homology are all the union of a surface with boundary in one half of the
cobordism with a surface with boundary in the other half.  The surface in the
first half of the cobordism is drawn in Figure~4.

\clearpage

\begin{figure}[h]
\vspace{2in}
\caption{The surface.} \label{fig4}
\end{figure}

The core of the two handle is also part of the surface.  For the other half of
the surface, note that the boundaries of the $2$-handles labeled $x,y,z$ and
$w$ are given by:
\[
\begin{array}{rll}
\p x &= &H + 2Q_1 \\
\p y &= &H + 3Q_2 \\
\p z &= &H + 5Q_3 \\
\p w &= &H + Q_1 + Q_2 + Q_3.
\end{array}
\]
After solving, we see that there is a surface with boundary the regular fiber,
given by
\[
15x + 10y + 6z - 30w.
\]
The surface has no self intersections in the standard half of the cobordism and
has $-30$ intersections in the $[0,1] \x \Sigma(2,3,5)$ part of the cobordism.
This implies that $\mbox{Sign } W = -1$.

When we compute the twisted signature, the surface in the standard half will be
the same as it was in the untwisted case because the representation is trivial
on that half.

The homology with twisted coefficients of a manifold, $W$, may be computed from
the universal cover, ${\tilde W}$.  The twisted homology is just the homology
of the complex:
\[
\to C_n({\tilde W}) \otimes_{{\Bbb Z}\pi_1(W)} su_2 \to C_{n-1}({\tilde W})
\otimes_{{\Bbb Z}\pi_1(W)} su_2 \to .
\]
The universal cover may be identified with a space of equivalence classes of
paths starting at a fixed base point.  A map from a cell to the universal cover
is determined by the induced map into $W$ together with one path from the base
point to the cell.  The group $\pi_1(W)$ acts on $C_*({\tilde W})$ by deck
transformations and it acts on $SU_2$ by the adjoint representation.  In
Figure~5, we have drawn a picture of a chain with boundary $j$ times the
regular fiber.

\clearpage

\begin{figure}[h]
\vspace{3in}
\caption{A twisted cocycle.} \label{fig5}
\end{figure}

The figures are slices $\{t\} \x \Sigma(2,3,5)$.  The surgery description of
$\Sigma(2,3,5)$ is drawn as a dashed set of curves.  The solid lines represent
slices of the cells representing the chain.  At the $t = 0$ slice, the $2$-cell
is a disk on the outside of the $2$-handle.  The solid curve is the boundary of
this disk.  It is oriented by the induced orientation.  Until $t = 1/3$ this
curve traces out an isotopy.  At $t = 1/3$ two more $2$-cells appear.  By
taking one pair of sides of each of the new rectangles, we get the slice at $t
= 1/6$.  By taking the other sides, we get the slice at $t = 1/2$.  It is
important to understand the labels on these cells.  Following the original path
from the base point might lead one to expect that the path to the top new
$2$-cell should wrap once around the meridian of the $2$-framed component.  At
first this is the case, but we change the path by an element of $\pi_1(W)$ and
act on the label by the inverse element.  In this case we compute the new label
by $\a(Q_1^{-1})j \a(Q_1) = -iji = -j$.  At $t = 2/3$ the final $2$-cell
appears, canceling the two meridians and leaving $j$ times the regular fiber.
If the coefficients were not twisted, the two meridians would have the wrong
orientations and would not cancel.  Drawing a parallel copy of this $2$-chain
shows that it has self intersection zero.  After drawing two more chains, we
can see that the twisted signature is trivial.

We will see that the fundamental representation extends over the cobordism,
thus Lemma~2 implies that
\[
\rho_M(\a) - \rho_{\Sigma(2,3,5)}(\a) = 3 \cdot \mbox{ Sign } W - \mbox{ Sign }
Q_{\a} = -3.
\]

The next step in the computation of the rho invariant is to construct a path of
representations on the new manifold.  The fundamental group of the new manifold
is:
\[
\pi_1(M) = \< H,Q_1,Q_2,Q_3,A_1,A_2 \mid [H,Q_i], A_kHA_k^{-1}H, Q_1^2H,
Q_2^3H, Q_3^5H, HQ_1Q_2Q_3A_1^2A_2^2\>.
\]
The group of the cobordism is:  $\pi_1(M)/\<A_1A_2\>$.  Define a path of
representations on $M$ by
\[
\begin{array}{ccc}
&H \mapsto -1, A_1 \mapsto \exp((t^2 - 1)\pi i/60), A_2 \mapsto 1, 
Q_1 \mapsto i, Q_2 \mapsto g_t \exp(\pi i/3)g_t^{-1}, \\ 
&Q_3 \mapsto g_t \exp(-\pi i/3)g_t^{-1} \exp(\pi i/2 + (1-t^2)\pi i/30).
\end{array}
\]
At $t = 1$ this representation extends over the cobordism and induces
the fundamental representation on $\Sigma(2,3,5)$.

To compute the group cohomology along this path of representations we must
compute the space of cocycles, $Z^1(M)$.  By the cocycle condition, we see that
$f(1) = 0$ and $f(g^{-1}) = -g^{-1}f(g)$ for any cocycle, $f$.  It follows that
a cocycle is determined by its values on a set of generators for $\pi_1(M)$.
Let $f(H) = h$ and so on.  With this conversion, we see that any cocycle must
satisfy the following equations, derived from the relations in the group.
\[
\begin{array}{rll}
0 &= &h + Hq_i - Q_ih - q_i = (1-Q_i)h, \\
0 &= &h + Ha_k + HA_kh - a_k = (1 + A_k)h, \\
0 &= &H(Q_3^4 + Q_3^3 + Q_3^2 + Q_3^1 + 1)q_3 + h = \frac {Q_3^5-1}{Q_3-1}
q_3 + h \\
0 &= &\frac {Q_2^3 - 1}{Q_2-1}q_2 + h \\
0 &= &\frac {Q_1^2 - 1}{Q_1-1}q_1 + h \\
0 &= &h + q_1 + Q_1q_2 + Q_1Q_2q_3 + Q_1Q_2Q_3(1 + A_1)a_1 + Q_1Q_2Q_3A_1^2(1 +
A_2)a_2.
\end{array}
\]
With the understanding that the lower case letters represent purely imaginary
quaternions and that multiplication by a group element just conjugates the
quaternion by the image of the group element in $SU_2 \cong Sp_1$, we see that
this is a linear system of $27$ equations and $18$ unknowns.  By solving, we
see that $\dim Z^1(M) = 9$ independent of $t$.  This shows that this path of
representations satisfies the hypothesis of Lemma~3 with $h = 6$.  This $h = 6$
represents the formal dimension of the moduli space of representations, so
there are many paths connecting the irreducible representation with a reducible
representation.  A similar computation shows that $\dim H^0(\Sigma(2,3,5);\a) =
\dim H^1(\Sigma(2,3,5);\a) = 0$.

It only remains to compute the rho invariant of the reducible representation.
The reducible representation, $\a_0: \pi_1(M) \to SU_2$ induces a $120 = 2^3
\cdot 3 \cdot 5$-fold cover of $M$ which is a circle bundle over a
non-orientable surface, ${\hat M} \to M$.  At the level of base orbifolds, the
cover has $2^2 \cdot 3 \cdot 5$ branch points of order $2$, $2^3\cdot 5$ branch
points of order $3$, and $2^3 \cdot 3$ branch points of order $5$.  A good
toy model of the desk transformations in this cover is given in Figure~6.

\clearpage

\begin{figure}[h]
\vspace{3in}
\caption{A branched cover.} \label{fig6}
\end{figure}

In this case, the action of the group of deck transformations extends over the
disk bundle, $W$, with boundary, $\p W = {\hat M}$.  Clearly $W$ is homotopy
equivalent to the non-orientable base orbifold.  Thus $H^2(W) = 0$, so
$H_{L^2}^2(W) = 0$, and $\mbox{Sign}(g,W) = 0$.  The image of $\a_0$ is
generated by $\exp(-\pi i/60)$.  The adjoint action is given by:
\[
\mbox{Ad}(e^{i\th}) \cdot i = i,\mbox{Ad}(e^{i\th})\cdot j = e^{2i\th}j,
\mbox{Ad}(e^{i\th}) \cdot k = e^{2i\th}k.
\]
Thus, $\mbox{Tr}(\mbox{Ad}(e^{i\th})) - 3 = 2(\cos 2\th - 1)$.

The fixed points of order $5$ come from $\exp(-\pi i/60)^{2^3\cdot 3\cdot n} =
\exp(-\frac {2\pi n}{5} i)$.  These elements correspond to sending $Q_3 \mapsto
\exp(\pi i/5)^{2^3\cdot 3 \cdot n} = \exp(\frac {4\pi n}{5} i)$.  From this we
see that the order $5$ fixed points have rotation angles; $\th_1 = \frac {4\pi
n}{5}$, $\th_2 = \frac {4\pi n}{5}$. It follows that the contribution to the
rho invariant from the fixed points of order $5$ is
\[
-\frac {2^3\cdot 3}{120} \sum_{n=1}^4 \cot(\frac {2\pi n}{5}) \cot(\frac {2\pi
n}{5} ) \cdot 2 \cdot (\cos(-\frac {4\pi n}{5}) - 1) = \frac {2^4 \cdot
3^2}{120} .
\]
In the same way, we can add up the contributions from all of the fixed points
to get:
\[
\rho_M(\a_0) = \frac {28}{15} .
\]
Putting this together with Lemma~2 and Lemma~3 gives, $\rho_{\Sigma(2,3,5)}(\a)
= \frac {28}{15} + 3 = \frac {73}{15}$.

In [A] we computed $CS_{\Sigma(2,3,5)} (\a) = \frac {1}{120}$.  Filling
in Lemma~1 gives the final answer:  $SF(\th,\a) = 1$.

In order to understand the meaning of this, extend the path of operators,
$B(A_t)$ to the whole real line by $B(A_t) = B(A_0)$ for $t < 0$ and $B(A_t) =
B(A_1)$ for $t > 1$.  Now solve the differential equation,
\[
\frac {du}{dt} + B(A_t) \cdot u = 0,
\]
by splitting $u$ into eigenfunctions.  The solution will decay (or grow) like
$e^{-\l_0t}$ as $t \to -\i$ and like $e^{-\l_1t}$ as $t \to \i$.  In other
words, the solution decays on the ends of ${\Bbb R}\x \Sigma(2,3,5)$ if and
only if $\l_0 < 0$ and $\l_1 > 0$.  The above differential equation is the
linearized version of the anti-self duality equations.  A moment's reflection
indicates that the spectral flow from $\th$ to $\a$ is the dimension of the
moduli space of anti-self dual connections on ${\Bbb R} \x \Sigma(2,3,5)$ with
boundary values $\th$ as $t \to -\i$ and $\a$ as $t \to \i$.  Starting with the
standard instanton on ${\Bbb R} \x S^3$, we can construct a one parameter
family of ASD connections on ${\Bbb R} \x \Sigma(2,3,5)$.  First, mod out by
the binary icosahedral group then translate that solution in the ${\Bbb R}$
direction.  David Austin has shown that these solutions make up the entire
moduli space of finite energy ASD connections on ${\Bbb R} \x \Sigma(2,3,5)$
with the given boundary conditions [Au].
\smallskip

We will now work through a hyperbolic example.  This manifold in Figure~7 is
hyperbolic according to the computer program snap pea [W].

\begin{figure}[h]
\vspace{2.6in}
\caption{The manifold, $M$.} \label{fig7}
\end{figure}

The fundamental group is generated by Wittinger generators,
\[
V_1,V_2,V_3,W_1,W_2,W_3,X_1,X_2,Y_1,Y_2,Z_1,Z_2
\]
with relations coming from
the crossings and the longitudes of the surgeries.  The crossing relation from
the upper left is, $V_3 = X_1V_1X_1^{-1}$.  Using the crossing relations, it is
possible to see that $\pi_1(M)$ is generated by $V_1,W_1,X_1,Y_1$, and $Z_1$.
The surgery relations may also easily be read from Figure~7.  For example the
surgery relation arising from the $V$ component is, $V^2ZYX = 1$.  It is an
elementary exercise to check that the assignments $V \mapsto -1$, $W \mapsto
-1$, $X \mapsto i \exp (-5/3 \pi k)$, $Y \mapsto i$, and $Z \mapsto \exp(2\pi
k/3)$ satisfy all of the crossing relations and all of the surgery relations.
The above assignment, therefore, defines a representation, $\a: \pi_1(M) \to
Sp_1$.  We will compute the spectral flow from the trivial representation to
$\a$.

The first step is to construct the flat cobordism, $W$, described in Lemma~2.
The new boundary component, $M_1$, is drawn in Figure~8.

\begin{figure}[h]
\vspace{3in}
\caption{A flat cobordant manifold, $M_1$.} \label{fig8}
\end{figure}

The fundamental group for $M_1$ may be obtained from the fundamental group of
$M$ by adding four new generators:  $A,B,C$, and $D$, four new relations:
$AVA^{-1}V$, $BVB^{-1}V$, $CWC^{-1}W$, $DWD^{-1}W$ and altering two of the
surgery relations to:  $V^2ZYXA^2B^2$ and
\[
W^4XV^{-1}XVX^{-1}Z^{-1}V^{-1}ZYZ^{-1}VC^2D^2.
\]
The induced representation on
$M_1$, $\a_1$, is clearly in the path component of a reducible representation,
$\a'_1$ on $M_1$.  To see this, just conjugate $\a_1(X)$ and $\a_1(Z)$ to $i$
and $\exp(2\pi i/3)$, set $\a'_1(B) = \a'_1(D) = 1$ and solve for $\a'_1(A)$
and $\a'_1(C)$ with the altered surgery relations.

The reducible representation is in the path component of the trivial
representation in the representation space of the link complement.  We may,
therefore, compute the Chern--Simons invariant with the Kirk and Klassen
formula [KK1], [A].  For each meridian pick a ``linear'' path
from the trivial element of $Sp_1$ and a positively oriented longitude.  As an
example, send $X$ to $\exp(2\pi ia_X(t))$ where $a_X(t) = \frac {1}{4} t$ and
notice that the longitude is:  $\l_X = X^4VW$, so that $b_X(t) = 4a_X(t) +
a_V(t) + a_W(t) = 2t$.  Plugging in all of the meridians and longitudes gives:
\[
CS_M(\a) \equiv -1/3\mod {\Bbb Z}.
\]

Elementary linear algebra will now be sufficient to compute the twisted
cohomology groups and the change in the rho invariant from $M,\a$ to
$M_1,\a'_1$.  By inspection, $\a$ is irreducible and $\th$ is trivial, so
\[
\dim H^0(M;\th) = 3 \mbox{ and } \dim H^0(M;\a) = 0.
\]
The group, $H^1(M;\th) \cong H_1(M;\th) \cong \pi_1(M)/[\pi_1(M),\pi_1(M)]
\otimes_{\Bbb Z} {\Bbb R}^3 \cong H_1(M;{\Bbb R}) \otimes {\Bbb R}^3$ as will
be verified when we compute $\mbox{Sign } W$.  According to
Lemma~\ref{lem3}, the
only contribution to $\mbox{Sign } W$ comes from $\mbox{coker}(H^1(M) \oplus
H^1(2V) \to H^1(2S^1 \x D^2))$, or by the universal coefficient theorem from
$\ker(H_1(2S^1 \x D^2) \to H_1(M) \oplus H_1(V))$.  The group, $H_1(2S^1 \x
D^2)$ is generated by $V$ and $W$.  Both $V$ and $W$ map to zero in $H_1(2V)$
as in Figure~4.  In $H_1(M)$, we have
\[
\begin{array}{rll}
\p a &= &4X + V + W \\
\p b &= &4Y + V + W \\
\p c &= &3Z + V - W \\
\p d &= &2V + X + Y + Z
\end{array}
\]
and
\[
\p e = 4W + X + Y - Z.
\]
From these equations we may solve for each of $V,W,X,Y$ and $Z$ in terms of
boundaries so that $H_1(M;{\Bbb R}) = 0$.  In particular, $V = \p(\frac {1}{66}
(57d - 18c + 3e - 15a - 15b))$, so we get a $2$-dimensional homology class
associated to $v$, namely $\frac {1}{66} (57d - 18c + 3e - 15a - 15b) - \frac
{1}{2} f$, where $f$ is the surface in Figure~4.  From this we may compute $V
\cdot V$ given $c \cdot c = 3$ etc.  After computing $W \cdot W$ and $V \cdot
W$, we see that
\[
\mbox{Sign } W = 2.
\]

A direct computation will show that
\[
\dim H^1(M;\a) = 3 \mbox{ and recall that } \dim H^1(M;\th) = 0.
\]
Furthermore, $\dim \mbox{ coker}(H^1(M;\a) \oplus H^1(V;\a) \to H^1(2S^1 \x
D^2;\a)) = 3$ and the coker is generated by $V - W$.  The corresponding labeled
surfaces in $W$ are a union of a pair of thin pants and a pair of fat pants
labeled with an $i$ around the $X$-component, also the same surface labeled
with either a $j$ or a $k$ around the $Y$-component.  See Figure~9.

\clearpage

\begin{figure}[h]
\vspace{2.3in}
\caption{Thin pants, fat pants.} \label{fig9}
\end{figure}

This shows that
\[
\mbox{Sign } Q_{\a} = 0.
\]
We will do a similar computation in a bit more detail later.

More linear algebra will verify that the path of representations from $\a'_1$
to $\a_1$ satisfies the hypothesis of Lemma~3.  Combining Lemma~2 and Lemma~3
gives:
\[
\rho_M(\a) = \rho_{M_1}(\a'_1) - 6.
\]
After computing $\rho_{M_1}(\a'_1)$, we will know all of the terms in the
formula for the spectral flow.  We could directly construct a $4$-manifold
which equivariantly bounds the cover of $M_1$ induced from $\a'_1$ and compute
$\rho_{M_1}(\a'_1)$ with Lemma~4, but it is easier to first perform a sequence
of flat cobordisms and a trivial deformation.

One natural way to build a flat cobordism is to attach a $2$-handle to $I \x M$
along a curve in the kernel of the representation.  Lemma~3 is also valid for
this type of cobordism, because the only fact about $V$ used in the proof was
that the second cohomology group with the appropriate coefficients is trivial.
For the first cobordism, $W_1$, just attach a $0$-framed two handle to the
meridian, $B$.  For $W_2$ attach a $0$-framed two handle to $D$.  Always call
the new boundary component of $W_n$, $M_{n+1}$.  The third cobordism is
constructed by sliding the $W$ handle over the $V$ handle in $M_3$ and then
attaching a $0$-framed two handle to the meridian of the $V$-handle.  See
Figure~10.

\clearpage

\begin{figure}[h]
\vspace{2.5in}
\caption{The cobordism $W_3$.} \label{fig10}
\end{figure}

In a surgery description, the linking matrix is the same as the intersection
matrix of the $4$-manifold.  Thus when we slide handle $W$ over handle $V$ the
matrix changes as though we replaced $W$ by $W + V$ in the basis of a quadratic
form.  The image under the representation of the handle that we slide does not
change.  However, the meridian to $V$ now represents $VW^{-1}$ in $\pi_1$.
Thus the loop labeled $U$ in Figure~10 is in the kernel of the representation.

Sliding the $Y$-handle over the $X$-handle and attaching a $0$-framed
$2$-handle to the meridian of $X$ produces the next cobordism, $W_4$.  See
Figure~11.

\begin{figure}[h]
\vspace{2.5in}
\caption{The cobordism $W_4$.} \label{fig11}
\end{figure}

At this stage we will deform the representation along the path $A \mapsto
\exp(2\pi i(\frac {1}{12} t - \frac {1}{2}))$, $C \mapsto \exp(2\pi i(-\frac
{1}{12} t))$ leaving all other meridians constant.  It is easy to see that all
of these representations are abelian but not central.  In fact, one can check
that the dimension of the group cohomology is constant along this path of
representations.  The proof of Lemma~4 shows that
\[
\rho_{M_5}(\a'_5) = \rho_{M_5}(\a_5).
\]

Continuing in this manner, we kill the $C$ component with the $W_5$ cobordism,
then slide the $A$ component over the $W$-component and kill the $W$-component
with the $W_6$ cobordism.  This leaves us with the manifold and representation
pictured in Figure~12.

\begin{figure}[h]
\vspace{2.5in}
\caption{The manifold, $M_7$.} \label{fig12}
\end{figure}

In order to apply Lemma~2 to compute the change in the rho invariant in this
sequence of flat cobordisms, we need to compute the signatures and twisted
signatures of the cobordisms.  The ordinary signature may be computed from the
linking matrix of the left end of the cobordism.  If the attaching map of the
new $2$-handle is trivial in $H_1(M;{\Bbb R})$, then the computation is exactly
the same as the computation of $\mbox{sign } W$ above.  If the attaching map is
non-trivial, then the signature is zero.  The $W_5$ cobordism is an example of
this second situation.  The linking matrix of the given surgery description of
$M_5$ is:
\[
\left[ \begin{array}{cccc}
0 & 2 & 0 & 0 \\
2 & 6 & 2 & 0 \\
0 & 2 & 0 & 0 \\
0 & 0 & 0 & 8
\end{array} \right]
\]
in the basis $\{A,W,C,Y\}$.  From this it follows that $H_1(M;{\Bbb R}) \cong
{\Bbb R}$ generated by $[C]$.  Thus, the kernel is trivial.  The result of
these computations is:
\[
\begin{array}{rll}
\mbox{Sign } W_1 = 0 &\mbox{Sign } W_2 = 0 &\mbox{Sign } W_3 = 0 \\
\mbox{Sign } W_4 = 1 &\mbox{Sign } W_5 = 0 &\mbox{Sign } W_6 = -1. 
\end{array}
\]

The twisted signatures may be computed in basically the same way.  Since the
coefficients are reducible, $H_1(M_k;\a) \cong H_1(M_k;{\Bbb R}) \oplus
H_1(M_k;{\tilde {\Bbb C}})$.  It follows that $\mbox{Sign}^Q W_k =
\mbox{ Sign } W_k + \mbox{ Sign}^{\Bbb C} W_k$.  The group $\ker(H_1(S^1 \x
D^2;{\tilde {\Bbb C}}) \to H_1(M_k;{\tilde {\Bbb C}}))$ is trivial for $k =
1,2,4$ and $5$, thus $\mbox{Sign}^Q W_k = \mbox{ Sign } W_k$ for these values
of $k$.  If $\a(z) \in S^1$ then $\a(z)$ acts trivially on the $i{\Bbb R}$ part
of $Sp_1$ and acts with weight two on the ${\Bbb C}j$ part, since
\[
e^{i\th} \cdot zj = e^{i\th}zje^{-i\th} = e^{2i\th}zj.
\]

The group, $H_1(M_3;{\tilde {\Bbb C}})$ is generated by the generators of
$\pi_1(M_3)$ with relations given by the Fox derivatives of the relators.  The
relation coming from the upper right crossing in Figure~10 is:  $UZU^{-1}Z^{-1}
= 1$.  Taking the Fox derivative gives:
\[
(1-e^{4\pi i/3})u = u + \a_3(U)^2z - \a(UZU^{-1})^2u - z = 0.
\]
A representative generating $H_2^{L^2}(W_3)$ is drawn in Figure~13.  It is
clear that it has trivial self intersection, so that $\mbox{Sign}^{\Bbb C} W_3
= 0$.

\begin{figure}[h]
\vspace{2in}
\caption{A representative for $U$.} \label{fig13}
\end{figure}

The signature, $\mbox{Sign}^{\Bbb C} W_6$ may be computed by the
method above.  The generator of $H_2^{L^2}(W_6;{\tilde {\Bbb C}})$ is fairly
complicated and the picture is messy.  There is an alternative way to compute
this signature.  First note that we are considering $H_2^{L^2}(W_6;{\tilde
{\Bbb C}})$ as a real vector space.  This means that $\mbox{Sign}^{\tilde {\Bbb
C}} W_6 = 2 \mbox{ Sign}_{\a} W_6$, where $\mbox{Sign}_{\a} W_6$ is the
signature, considering $H_2^{L^2}(W_6;{\tilde {\Bbb C}})$ as a complex vector
space.  Inventing formula $2.5$ from [APS] gives:
\[
\mbox{Sign}_{\a} X = \frac {1}{|G|} \sum_{g \in G} \mbox{ Sign}(g,{\tilde X})
\chi_{\a}(g).
\]
In our situation, $G = {\Bbb Z}_2$, $X = W_6$ and ${\tilde X}$ is the $2$-fold
cover of $W_6$ induced from the representation.  By restricting the
intersection form to the positive and negative eigenspaces of the generator of
${\Bbb Z}_2$, we get:
\[
\mbox{Sign}(1,{\tilde X}) = \mbox{ Sign}({\tilde X}) = \mbox{ Sign } Q|_{E_+} +
\mbox{ Sign } Q|_{E_-}
\]
and
\[
\mbox{Sign}(-1,{\tilde X}) = \mbox{ Sign } Q|_{E_+} - \mbox{ Sign } Q|_{E_-}.
\]
Furthermore, a transfer argument shows that $\mbox{Sign } Q|_{E_+} =
\mbox{ Sign } W_6$, giving:
\[
\mbox{Sign}_{\a} W_6 = \mbox{ Sign } {\tilde X} - \mbox{ Sign } W_6 = -2 - (-1)
= -1.
\]
The cover, ${\tilde X}$ is pictured in Figure~14.

\begin{figure}[h]
\vspace{3in}
\caption{The cover ${\tilde X}$.} \label{fig14}
\end{figure}

Summarizing:
\[
\begin{array}{rll}
\mbox{Sign}^Q W_1 = 0 &\mbox{Sign}^Q W_2 = 0 &\mbox{Sign}^Q W_3 = 0 \\
\mbox{Sign}^Q W_4 = 1 &\mbox{Sign}^Q W_5 = 0 &\mbox{Sign}^Q W_6 = -3
\end{array}
$$
so that
\[
\rho_{M_7}(\a'_1) = \rho_{M_7}(\a_7) - 2.
\]
by Lemma~2.

We are now in a position to finish the computation.  Since $M_7 = M_8 \#
L(3,-1)$,
\[
\rho_{M_7}(\a_7) = \rho_{M_8}(\a_7|_{M_8}) + \rho_{L(3,-1)}(\a_7|_{L(3,-1)}).
\]
The induced cover of the lens space is $S^3$ which equivariantly bounds $D^4$
with one fixed point.  Lemma~4, therefore, implies that
$\rho_{L(3,-1)}(\a_7|_{L(3,-1)}) = -2/3$.  In order to compute
$\rho_{M_8}(\a_7|_{M_8})$, we need to construct a $4$-manifold with a ${\Bbb
Z}_4$ action so that the quotient map on the boundary is exactly the covering
of $M_8$ induced from the representation.

There are two general methods to construct such a $4$-manifold.  The first
method would be to use Kirby moves to produce a surgery description of $M_8$ so
that the meridian of one unknotted component is sent to a generator and the
meridians of all other components are in the kernel of the representation.  The
unknotted component alone would produce a lens space.  The cover of this lens
space is $S^3$ which bounds $D^4$.  By gluing $2$-handles to the lifts of the
remaining components, we will produce the desired $4$-manifold.  The second
method is described in the work of Casson and Gordon [CG].  It is this
second method that we will use here.  Start by performing handle slides until
all meridians are sent to the generator.  End the construction by taking the
branched cover of the surgery description branched along a pushed in Seifert
surface glued to the cores of all of the $2$-handles.

One apparent difference between our computations and the Casson--Gordon
computations is the difference between using $L^2$-cohomology or ordinary
homology to compute signatures.  There really is no difference.  By duality, we
could compute in $L^2$-homology, $H_2^{L^2}(W) \equiv \mbox{ Im}(H_2(W) \to
H_2(W,\p W))$.  If $F_1$ is any class in $H_2(W)$ which is trivial in
$H_2^{L^2}(W)$, then $F_1$ is represented by a surface in $\p W$.  Using a
collar of $\p W$, any other surface may be pushed into the interior of $W$.
Thus the intersection number with $F$ is zero and $F$ is in the null space of
the intersection form on $H_2(W)$.

Recapping from Casson and Gordon, if $\<\tau \mid \tau^m = 1\> = {\Bbb Z}_m$
acts on $X^4$, then
\[
\begin{array}{rll}
\mbox{Sign}(\tau^s,X) &= &\sum_{r=0}^{m-1} \exp(2\pi
irs/m)\mbox{Sign}(Q|_{E_r}) \\
&= &\mbox{Sign}(X/{\Bbb Z}_m) + \sum_{r=1}^{m-1} \exp(2\pi
irs/m)\mbox{Sign}(Q|_{E_r}).
\end{array}
\]
Here $E_r$ is the $\exp(2\pi ir/m)$-eigenspace of $\tau$ and the second line
is true by a transfer argument.  The transfer argument is just that $[\a] \in
E_0$ implies that $[\a] = [\frac {1}{m} \sum_{r=0}^{m-1} (\tau^r)^*\a]$, but
$1/m \sum_{r=0}^{m-1} (\tau^r)^*\a$ is equivariant and may therefore be pushed
forward.  This constructs an isomorphism, $E_0(X) \to H(X/{\Bbb Z}_m)$.  Casson
and Gordon use a Mayer--Vietoris argument to show that the signatures of the
other eigenspaces may be computed from the cyclic cover of the complement of a
pushed in Seifert surface.  If $\{x_k\}$ is a basis for the first homology of
the Seifert surface, then $\{\sum_{s=0}^{m-1} \exp(2\pi irs/m)\tau_*^sx_k^*\}$
is a basis for $E_r$.  In this basis, $Q|_{E_r} = (1 - \exp(-2\pi ir/m))A + (1
- \exp(2\pi ir/m))A^*$, where $A$ is the Seifert form.

Figure~15 shows the manifold, $M_8$ after the appropriate handle slides have
been performed.

\newpage

\begin{figure}[h]
\vspace{2in}
\caption{The manifold, $M_8$.} \label{fig15}
\end{figure}

Figures~16 and 17 show two parts of the Seifert surface so that the Seifert
form may be computed.

\begin{figure}[h]
\vspace{3in}
\caption{Top half of Seifert surface.} \label{fig16}
\end{figure}

\clearpage

\begin{figure}[h]
\vspace{3in}
\caption{Bottom half of Seifert surface.} \label{fig17}
\end{figure}

\clearpage
${}$


\noindent
Plugging into the preceding formulae for $g$-signatures, and then Lemma~3,
gives:
\[
\rho_{M_8}(\a_7|_{M_8}) = 150.
\]
To compute the signatures, we just entered
``$N[\mbox{Eigenvalues}[Q_{E_r}]]$'' into
{\em Mathematica} and counted the number of positive and negative eigenvalues
[Wo].  There are in fact very good numerical methods for computing
signatures.
See [?] for example.

Combining all of the previous computations with Lemma~1 completes the
computation to
give:
\[
SF_M(\th,\a) = 65.
\]

To conclude, we have explained a method to compute the spectral flow.  This
method
works for a large class of representations on a large class of manifolds.
Namely, it
works for any representation that is flat cobordant-flat deformation
equivalent to a
reducible representation.  It is conceivable that this method works for any
representation on any $3$-manifold.  The recent work of Reznikov gives some
slight evidence
supporting this possibility.

I would like to thank N.~Thurston for help with snap pea, P.~DeSouza for
help with
mathematica, and A.~Casson, R.~Fintushel, P.~Kirk and E.~Klassen for various
conversations.

\bigskip
\noindent
Department of Mathematics \\
University of California, Berkeley \\
Berkeley, CA 94720 \\
email: dav@@math.berkeley.edu


\begin{thebibliography}{KKZ}
\bibitem[APS]{AtPaSi} M. Atiyah, V. Patodi and I. Singer, {\em Spectral asymmetry
and Riemannian geometry I, II, III}, Math. Proc. Camb. Phil. Soc. {\bf 77,78,79}
(1975), 43--69, 405--432, 71--99.

\bibitem[A]{Auc} D. Auckly, {\em Topological methods to compute Chern--Simons
invariants}, Math. Proc. Camb. Phil. Soc. {\bf 115} (1994), 229--251.

\bibitem[Au]{Aus} D. Austin, {\em Equivariant Floer groups for binary polyhedral spaces}, Math. Ann. {\bf 302} (1995) no 2, 295--322.

\bibitem[CG]{CaGo} A. Casson and C. Gordon, {\em On Slice Knots in dimension
three}, Proc. Symp. Pure Math. {\bf 32} (1978), 39--53.

\bibitem[FS1]{FinSt1} R. Fintushel and R. Stern, {\em Pseudofree orbifolds},
Ann. of Math. {\bf 122} (1985), 335--364.

\bibitem[FS2]{FinSt2} R. Fintushel and R. Stern, {\em Instanton homology of
Seifert-fibered homology spheres}, Proc. Lond. Math. Soc. (3) {\bf 61} (1990),
109--137.

\bibitem[KK1]{KiKl1} P. Kirk and E. Klassen, {\em Chern--Simons invariants of
$3$-manifolds and representation spaces of knot groups}, Math. Ann. {\bf 287}
(1990), 347--367.

\bibitem[KK2]{KiKl2} P. Kirk and E. Klassen, {\em Computing spectral flow via
cup products}, J. Diff. Geom. {\bf 40} (1994), 505--562.

\bibitem[KKR]{KiKlRu} P. Kirk, E. Klassen and D. Ruberman, {\em Splitting the
spectral flow and the Alexander matrix}, Comment. Math. Helv. {\bf 69} (1994) 
no 3, 375--416.

\bibitem[R]{Re} A. Rezinikov, {\em Rationality of secondary classes}, J. Diff.
Geom. to appear.

\bibitem[W]{We} J. Weeks, ``Snap pea'', available at jweeks@@midd.bitnet.

\bibitem[Wo]{Wol} S. Wolfman, Mathematica, A system for doing mathematics by
computer, Addison--Wesley, 1988.

\end{thebibliography}
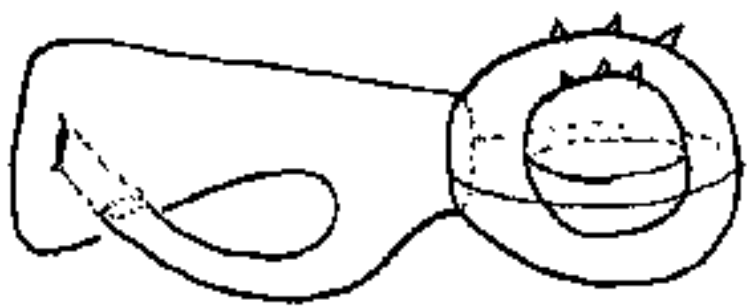
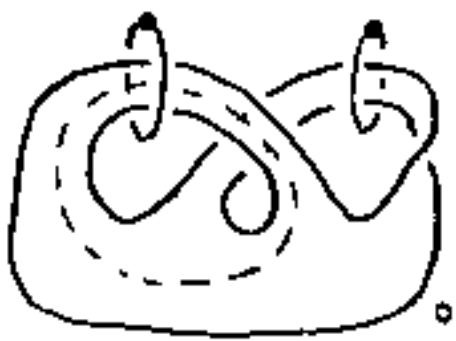
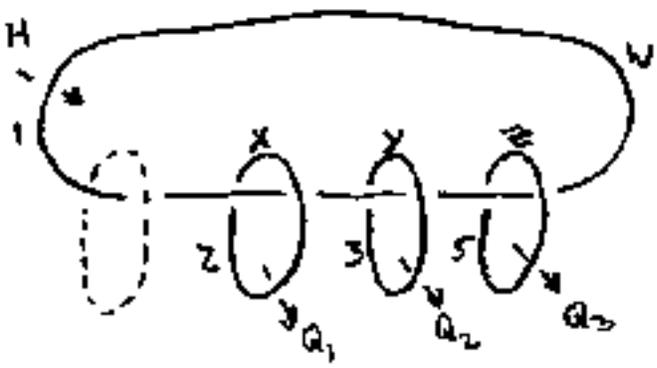
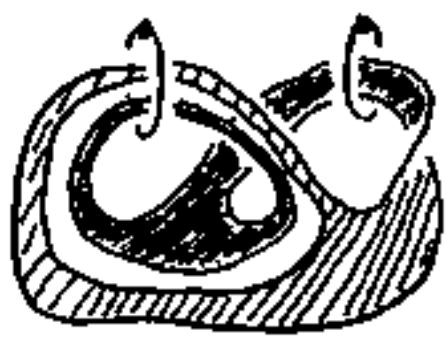
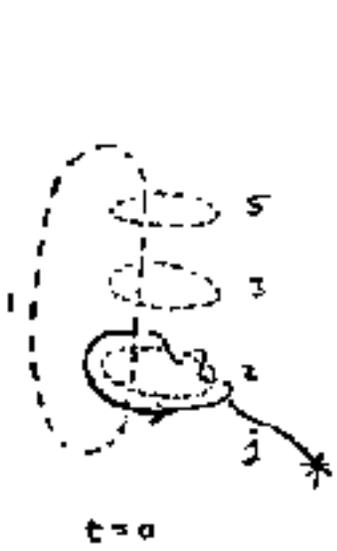
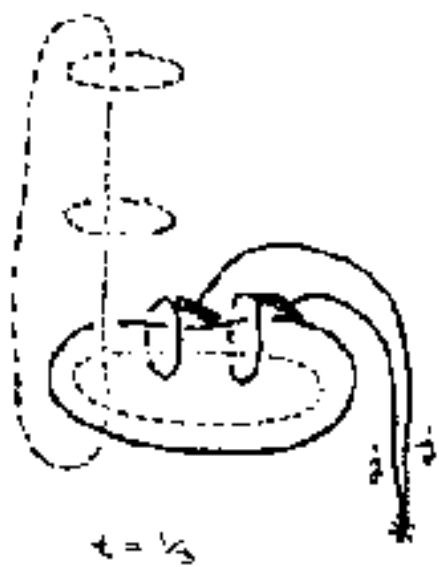
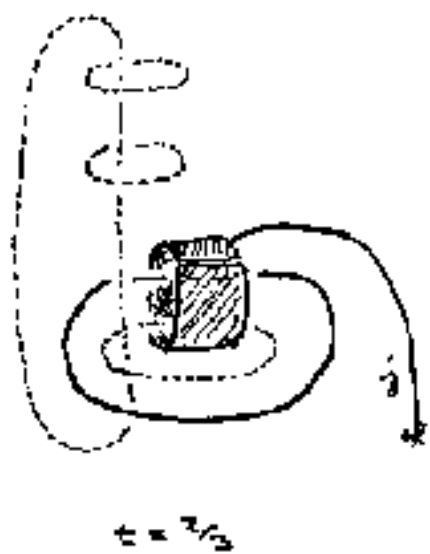
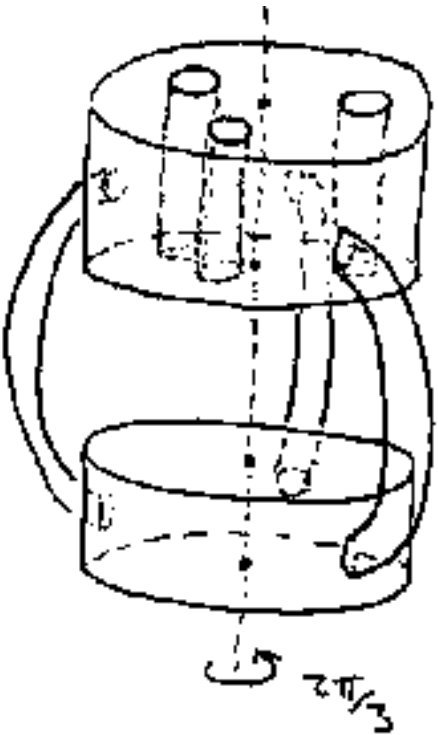
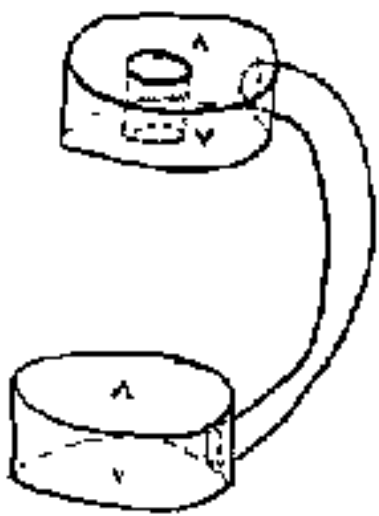
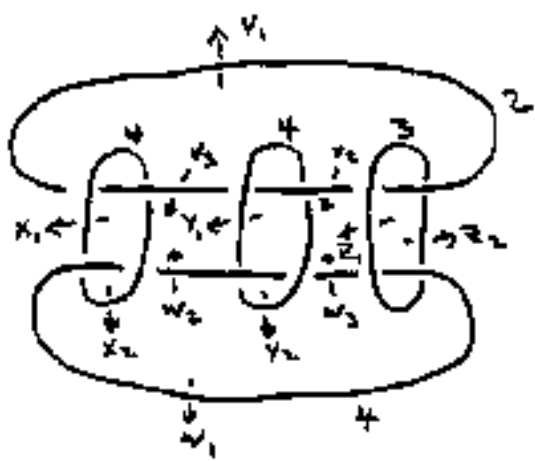
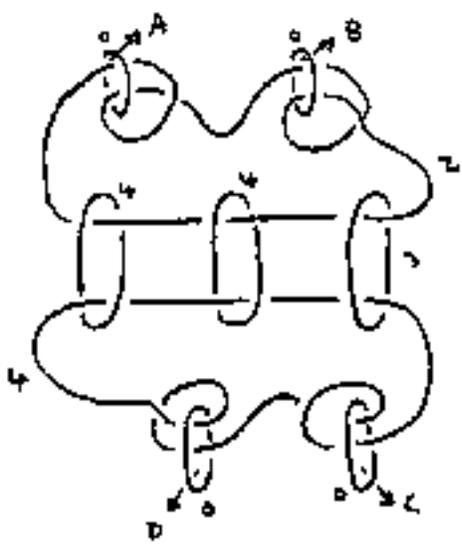
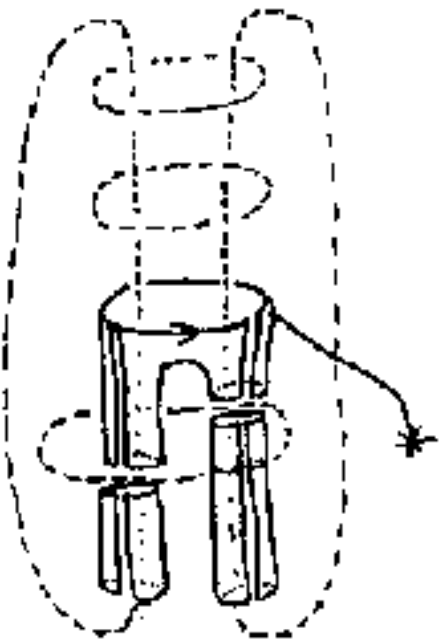
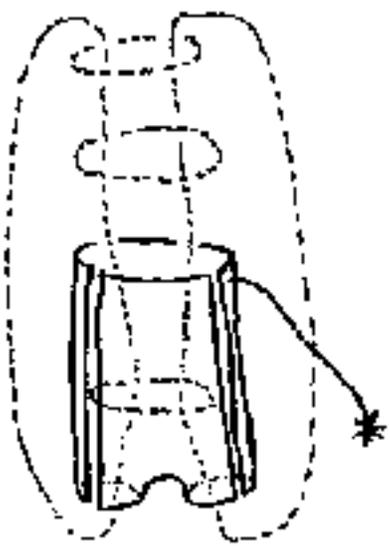
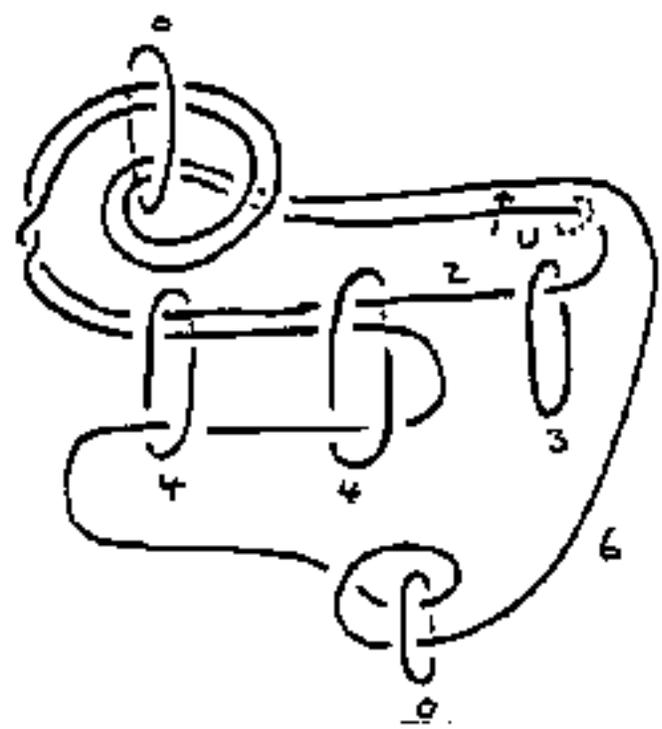
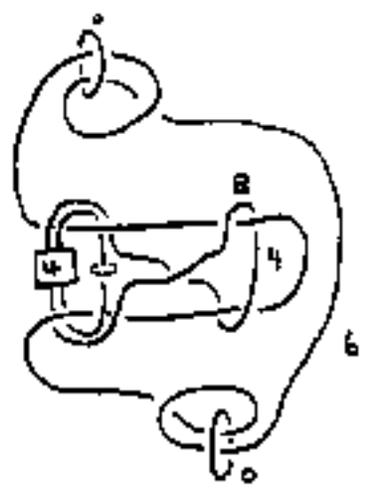
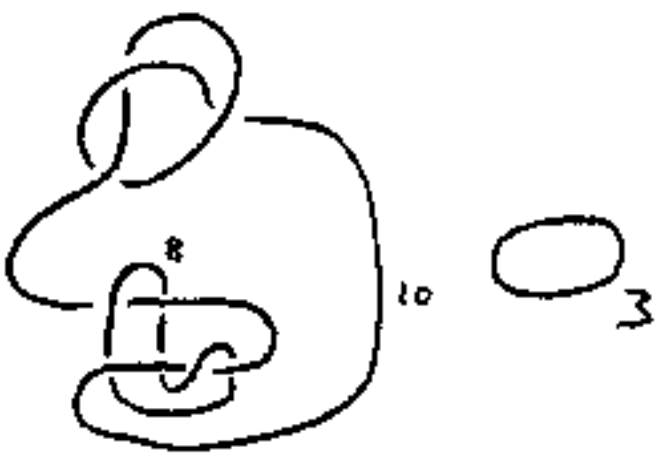
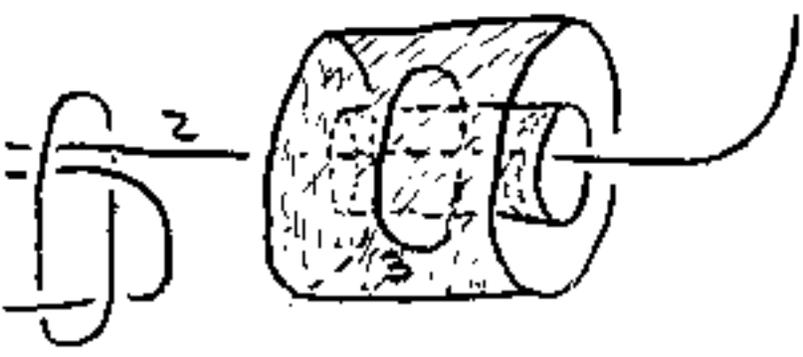
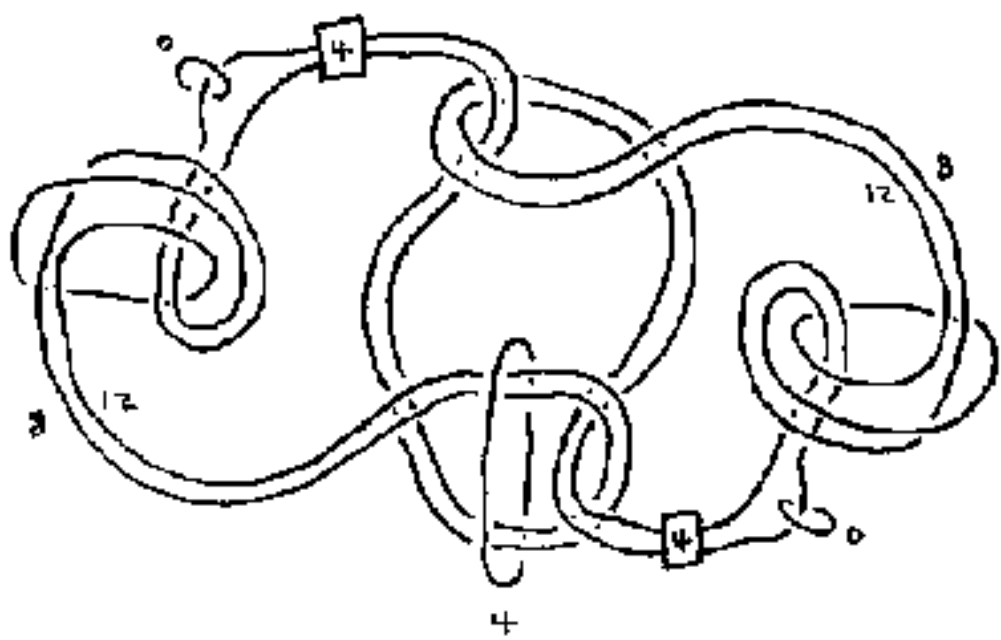
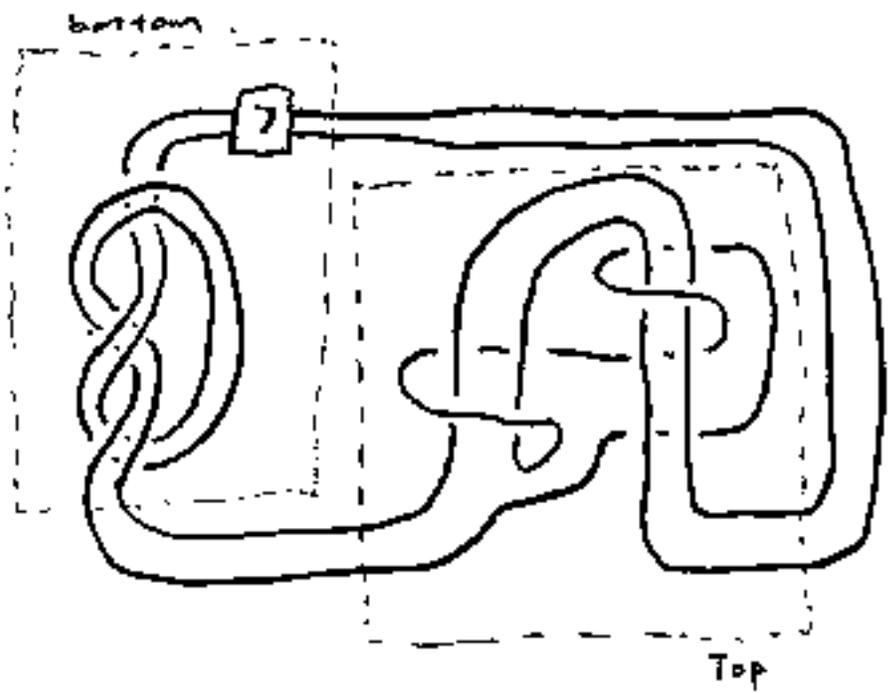
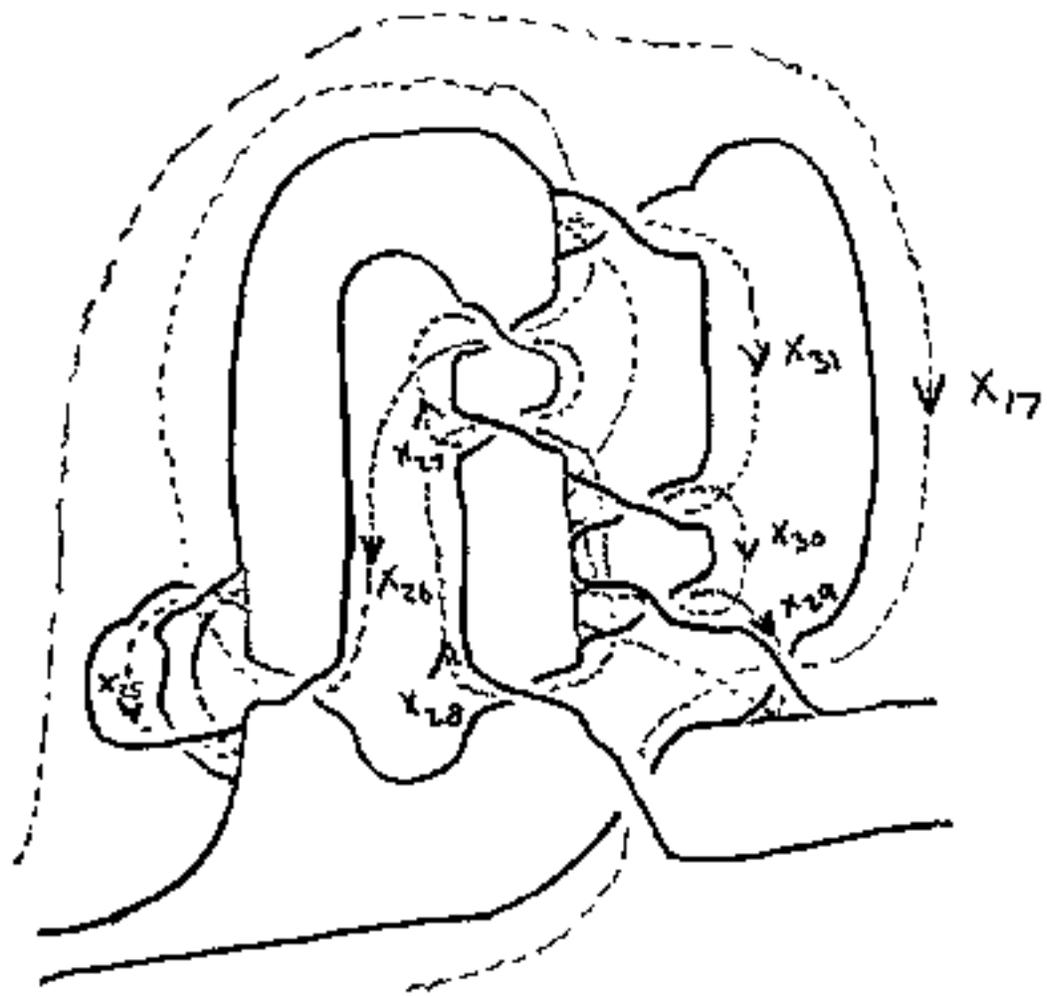
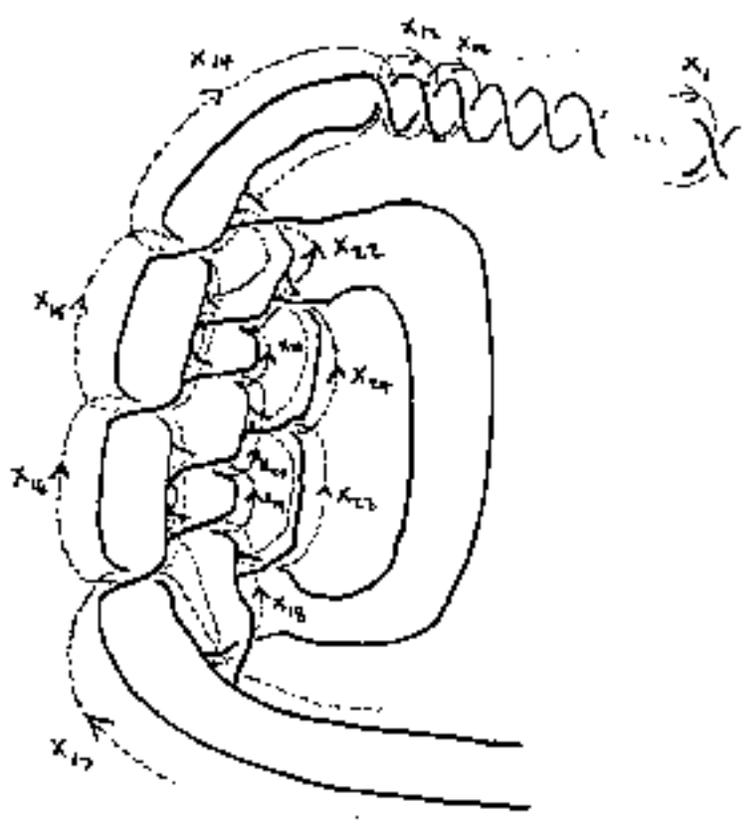
\end{document}